\documentclass[12pt]{article}
\usepackage{amsmath}
\usepackage{amssymb}
\usepackage{euscript}
\usepackage[latin1]{inputenc}

\normalbaselines
\newcommand{\qed}{\sqcap\kern-8pt\sqcup}

\newenvironment{proof}
      {\par\noindent{\it Proof\/: }\nopagebreak\normalsize}%
                                  {\linebreak[2]\hspace*{\fill}$\qed$\ifdim\lastskip<12pt
       \removelastskip \penalty-200  \vskip12pt  \fi}

\font\frten=eufm10 at 12pt
\font\freight=eufm10
\font\frsix=eufm8
\newfam\frfam\textfont\frfam=\frten
\scriptfont\frfam=\freight
\scriptscriptfont\frfam=\frsix

\newcommand{\CC}{\mathbb{C}}
\newcommand{\PP}{\mathbb{P}}
\newcommand{\RR}{\mathbb{R}}
\newcommand{\NN}{\mathbb{N}}
\newcommand{\ZZ}{\mathbb{Z}}

\newtheorem{thm}{Theorem}[section]
\newtheorem{prop}[thm]{Proposition}
\newtheorem{cor}[thm]{Corollary}
\newtheorem{lem}[thm]{Lemma}

\newtheorem{rema}[thm]{Remark}

\newtheorem{exem}[thm]{Examples}

\def \Clif{{\rm Clif}}
\def \dim{{\rm dim}}
\def \Cliff{{\rm Cliff}}
\def \O{{\cal O}}

\def \L{{\cal O}}

\def \air{{\vskip 12pt\noindent}}

\def \Div{{\rm Div}}
\def \Gal{{\rm Gal}}

\def \deg{{\rm deg}}

\begin{document}

\title{Clifford Theorem for real algebraic curves
\footnote{Work supported by the European Community's 
Human Potential Programme 
under contract HPRN-CT-2001-00271, RAAG.}}

\author{Jean-Philippe Monnier\\
       {\small D\'epartement de Mathématiques, Universit\'e d'Angers,}\\
{\small 2, Bd. Lavoisier, 49045 Angers cedex 01, France}\\
{\small e-mail: monnier@tonton.univ-angers.fr}}
\date{}
\maketitle

\begin{abstract}
{ We establish for smooth projective real curves the equivalent of the 
classical Clifford inequality known for complex curves. We also study
the cases when equality holds.}
\end{abstract}

\section*{Introduction}

In this note, a real algebraic curve $X$ is a smooth proper geometrically
integral scheme over $\RR$ of dimension $1$.
A closed point $P$ of $X$ 
will be called a real point if the residue field at $P$ is $\RR$, and 
a non-real point if the residue field at $P$ is $\CC$.
The set of real points, $X(\RR)$,
will always be assumed to be non empty.
It decomposes into 
finitely many connected components, whose number will be denoted by $s$.
By Harnack's Theorem we know that $s\leq g+1$,
where $g$ is the genus of $X$. 
A curve with $g+1-k$ real connected components
is called an $(M-k)$-curve. 

The group $\Div(X)$ of divisors on $X$ 
is the free abelian group generated by the 
closed points of $X$. If $D$ is a divisor on $X$, we will denote by
$\L (D)$
its associated invertible sheaf. The dimension of the space 
of global sections of this sheaf will be denoted
by $\ell (D)$. Let $D\in \Div(X)$, 
since a principal divisor has an even degree on each connected 
component of $X(\RR)$ (\cite{G-H} Lem. 4.1),
the number $\delta(D)$ (resp. $\beta(D)$) of connected components $C$
of $X(\RR )$
such that the degree of the restriction of $D$ to $C$ is odd (resp even), 
is an invariant of the linear system $|D|$
associated to $D$. Let $K$ be the canonical divisor. If 
$\ell (K-D)=\dim\, H^1 (X,\L (D) )>0$, $D$ is said to be special. If not,
$D$ is said to be non-special. By Riemann-Roch, if $\deg(D)>2g-2$ then $D$
is non-special. Assume $D$ is effective and let $d$ be its degree. If
$D$ is non-special then the dimension of the linear system $|D|$ is given by 
Riemann-Roch. If $D$ is special, then the dimension of the linear system
$|D|$ satisfies 
$$\dim\, |D|\leq \frac{1}{2} d.$$
This is the well known Clifford inequality for complex curves that
stands obviously for real curves. 

If $X$ is an $M$-curve or an $(M-1)$-curve, then Huisman (\cite[Th. 3.2]{Hu})
has shown that 
$$\dim\, |D|\leq \frac{1}{2} (d-\delta (D)).$$
We have already proved that Huisman's inequality
is also valid for almost all real hyperelliptic curves,
for example when $s\not= 2$. But 
there is a family of real hyperelliptic curves 
with $2$ real connected components for which 
there exist some special divisors $D$ satisfying
$\dim\, |D|= \frac{1}{2} d >\frac{1}{2} (d-\delta (D))$
(\cite[Th. 4.4, Th. 4.5]{Mo}).

In this note we establish a new Clifford inequality for real curves,
completing the inequality given by Huisman.

\air
{\bf Theorem A}
{\it Assume $D$ is effective and special. Then, either
$$\dim\, |D|\leq \frac{1}{2} (d-\delta (D))\,\,\,\, (\Clif \,1)$$
or 
$$\dim\, |D|\leq \frac{1}{2} (d-\beta (D))\,\,\,\, (\Clif \,2).$$ }
\air

The cases for which equality holds in (Clif 1) or (Clif 2) are also studied.
Moreover, looking at divisors that do not satisfy the inequality (Clif 1), we
obtain the following theorem.

\air
{\bf Theorem B}
{\it Let $X$ be a real curve.
Let $D\in \Div(X)$ be an effective and special divisor of degree $d$.
\begin{description}
\item[({\cal{i}})] If $X$ is hyperelliptic then either
$$\dim\, |D|\leq \frac{1}{2} (d-\delta (D))\,\,\,\, (\Clif \, 1)$$
or 
$$\dim\, |D|\leq \frac{1}{2} (d-\frac{1}{2}(s-2)),$$
\item[({\cal{ii}})] else, either
$$\dim\, |D|\leq \frac{1}{2} (d-\delta (D))\,\,\,\, (\Clif \, 1)$$
or
$$\dim\, |D|\leq \frac{1}{2} (d-\frac{1}{2}(s-1)).$$
\end{description}}
\air

The author wishes to express his thanks to D. Naie for several helpful
comments concerning the paper.

\section{ Preliminaries}

We recall here some classical concepts and more notation that we will
be using throughout the paper.

Let $X$ be a real curve. We will denote by $X_{\CC}$ the base extension of $X$
to $\CC$. 
The group $\Div(X_{\CC})$ of divisors on $X_{\CC}$
is the free abelian group on the 
closed points of $X_{\CC}$. The Galois group $\Gal(\CC /\RR)$
acts on the complex variety $X_{\CC}$ and also on $\Div(X_{\CC})$. 
We will always indicate this action by a bar.
If $P$ is a non-real point of $X$, identifying $\Div(X)$
and $\Div(X_{\CC})^{\Gal(\CC /\RR)}$, then
$P=Q+\bar{Q}$ 
with $Q$ a closed point of $X_{\CC}$. 
If $D$ is a divisor on $X_{\CC}$, we will denote by 
$\L (D)$
its associated invertible sheaf and by $\ell_{\CC} (D)$ 
the dimension of the space 
of global sections of this sheaf. If $D\in \Div(X)$,
then
$\ell (D) =\ell_{\CC} (D)$.

As usual, a $g_d^r$ is an $r$-dimensional complete linear system
of degree $d$ on $X$ (or $X_{\CC}$). A $g_d^r =|D|$ on $X$ is called 
base point free if $\ell_{\CC} (D-P)=\ell_{\CC} (D)-1=r$ for any closed point
$P$ of $X_{\CC}$. Then it defines a morphism $\varphi :\, X\rightarrow
\PP_{\RR}^r$ onto a non-degenerate (but maybe singular) curve in
$\PP_{\RR}^r$. If $\varphi$ is birational (resp. an isomorphism) onto
$\varphi (X)$, the $g_d^r$ (or $D$) is called simple (resp. very
ample). Let $X'$ be the normalization of $f(X)$, and assume $D$ is not
simple i.e. $D-P$ has a base point for any closed point $P$ of $X_{\CC}$. Thus,
the induced morphism $\varphi :X\rightarrow X'$ is a non-trivial covering
map of degree $k\geq 2$. In particular, there is $D'\in \Div(X')$ such that
$|D'|$ is a $g_{\frac{d}{k}}^{r}$ and such that $D=\varphi^* (D')$, i.e.
$D$ is induced by $X'$. If $g'$ denote the genus of $X'$,
$|D|$ is classically called compounded of an involution of order $k$
and genus $g'$. In the case $g'>0$, we speak of an irrational involution 
on $X$.

Let $D\in \Div(X)$ be an effective divisor with base points. Let $D'\in
\Div(X)$ be 
a subdivisor of $D$. We say that $D'$ is a base point free part
of $D$ if $D'$ is base point free and the degree of $D'$
is maximal with this property. If $D'$ is a base point free
part of $D$, we may write $D=D'+E$ with $E$ an effective divisor.
Let $k$ be the number of non-real points (counted with multiplicity)
in the support of $E$, then $$\dim\, |D|\leq\, \dim\, |D' | +k.$$

The reader is referred to \cite{ACGH} and
\cite{Ha} for more details on special divisors. Concerning real curve,
the reader may consult \cite{BCR} and \cite{G-H}.
For $0\leq a\in\RR$ we denote by $[a]$ the integer part of $a$, i.e. 
the biggest integer $\leq a$.

\section{Clifford inequalities for real curves}

Before proving Theorem A stated in the introduction, we need to
establish some
preliminary results.

\begin{lem}
\label{1parcomp}
Let $P_1 ,\ldots ,P_{\delta} $ be real points of $X$ 
such that no two of them belong to the same connected component
of $X(\RR )$. Let $D$ be the divisor $P_1 +\ldots + P_{\delta} $.
Then
\begin{description}
\item[1)] $\dim\, |D| =0$ if $\delta <s$ and,
\item[2)] $\dim\, |D| \leq 1$ if $\delta =s$, and in case $\dim\, |D| =
  1$ then $D$ is base point free.
\end{description}
\end{lem}

\begin{proof}
Assume $\delta <s$ and $\dim\, |D| >0$. Choose a real point $P$ in one
of the $s-\delta $ real connected components that do not contain 
any of the points $P_1 ,\ldots ,P_{\delta} $. Then $\L (D-P)$ has a global
section and $D-P$ should be linearly equivalent to an effective divisor
$D'$ of degree $\delta -1$ satisfying $ \delta (D') =\delta +1$.
This is impossible, proving 1).

Assume $\delta =s$, then $\dim\, |D|\leq 1+\dim\, |P_2 +\ldots + P_s |\leq 1$
by 1). Suppose $\dim\, |D| =1$ and $D$ is not base point free. If
$D$ has a real base point $P$, then $\dim\, |D-P|=1$ and 
$\deg (D-P)=\delta (D-P)=s-1$, contradicting 1). If $D$ has a non-real
base point $Q$, then $\ell (D-Q)>0$ and $D-Q$ 
should be linearly equivalent to an effective divisor
$D'$ of degree $s -2$ satisfying $ \delta (D') =s$,
which is again impossible.
\end{proof}

The following lemma is due to Huisman \cite{Hu}.
\begin{lem} 
\label{huisman}
Let $D\in \Div(X)$ be an effective divisor of degree $d$ and 
assume $d+\delta (D)< 2s$.
Then $$\dim\, |D|\leq \,\frac{1}{2} (d-\delta (D)).$$
\end{lem}

\begin{proof}
Set $k=\frac{1}{2} (d-\delta (D))$. Then $k\geq 0$, since $D$
is effective. We have $\delta (D) +k<s$ by the hypotheses.
Choose $P_1, \ldots , P_k $ real points among
the $\beta (D)$ real connected components on which the degree
of the restriction of $D$ is even, such that no two of these points
belong to the same real connected component. Let $D' =D- P_1 - \ldots - P_k$.
Then $\deg (D')=d-k = \delta (D')= \delta (D) +k<s$. 
By Lemma \ref{1parcomp}, if $\ell (D')>0$, then $\dim\, |D'|=0$. Finally,
$\dim\, |D|\leq \dim\, |D'|+k \leq \,\frac{1}{2} (d-\delta (D)).$
\end{proof}

Generalizing the previous lemma, we get:
\begin{lem}
\label{huisman2}
Let $D$ be an effective divisor of degree $d$ on a
real curve $X$. Assume that $d+\delta(D)\leq 2s +2k $ with
$k\in \NN$. Then $$\dim\, |D|\leq \frac{1}{2} (d-\delta (D))+k+1.$$
\end{lem}

\begin{proof}
We proceed by induction on $k$. We explain the case $k=0$, the general
induction step being completely similar. So, assume that
$d+\delta(D)=2s$.

If $\delta(D)=0$ then choose a non-real point $Q$. Since 
$\deg(D-Q)+\delta(D-Q)<2s$ and $\delta(D-Q)=\delta (D)=0$, if $\ell
(D-Q)>0$ then, using Lemma \ref{huisman}, we have 
$\dim\, |D-Q|\leq \frac{1}{2} (\deg(D-Q)-\delta (D-Q))=\frac{1}{2}
(d-2-\delta (D))$. 
Hence $\dim\, |D|\leq \frac{1}{2} (d-\delta (D))+1$.
If $\ell (D-Q)=0$ then $\dim\, |D|\leq 1 \leq \frac{1}{2} (2s)=\frac{1}{2} d
= \frac{1}{2} (d-\delta (D))$ since we have assumed $s$ to be $>0$.

If $\delta(D)\not =0$, then choose a real point $P$ in a real connected
component $C$ such that the degree of the restriction of $D$ to $C$
is odd. Since $\deg(D-P)+\delta(D-P)<2s$, if $\ell (D-P)>0$ then,
using Lemma \ref{huisman}, we get 
$\dim\, |D-P|\leq \frac{1}{2} (\deg(D-P)-\delta (D-P))=\frac{1}{2}
(d-\delta (D))$. 
Hence $\dim\, |D|\leq \frac{1}{2} (d-\delta (D))+1$. 
If $\ell (D-P)=0$ then $\dim\, |D|=0 \leq 
\frac{1}{2} (d-\delta (D))$ ($d\geq \delta (D)$ since $D$ is effective).
\end{proof}

The following lemma will allow us to restrict the study to base point free
linear systems.
\begin{lem} 
\label{reductbasepointfree}
Let $D\in \Div(X)$ be an effective divisor of degree $d$.
Let $D'$ be a degree $d'$ base point free part of $D$. 
\begin{description}
\item[({\cal{i}})] If
$$\dim\, |D' |\leq \frac{1}{2} (d'-\delta (D'))+k$$
for a positive integer $k$, then
$$\dim\, |D |\leq \frac{1}{2} (d-\delta(D))+k .$$
\item[({\cal{ii}})]
If $$\dim\, |D' |\leq \frac{1}{2} (d'-\beta (D'))-k$$
for a positive integer $k$, then
$$\dim\, |D |\leq \frac{1}{2} (d-\beta (D))-k .$$
\end{description}
\end{lem}

\begin{proof}
Assume $D'\in \Div(X)$ is an effective divisor of degree $d'$
satisfying $$\dim\, |D' |\leq \frac{1}{2} (d'-\delta (D'))+k$$
for a positive integer $k$.
Let $P$ be a real point such that $\dim\, |D' +P | =\dim\, |D'|$.
Then $\dim\, |D' +P | \leq \frac{1}{2} (d'-\delta (D'))+k \leq 
\frac{1}{2} (\deg(D'+P) -\delta (D'+P))+k $.
Let $Q$ be a non-real point such that $\dim\, |D' +Q | \leq \dim\, |D'|+1$.
Then $\dim\, |D' +Q | \leq \frac{1}{2} (d'-\delta (D'))+1+k =
\frac{1}{2} (\deg(D'+Q) -\delta (D'+Q))+k $.

Statement ({\cal{i}}) is a consequence of the above results.
For statement ({\cal{ii}}), the proof is similar.
\end{proof}

Let $D$ be a special divisor.
Recall that $\delta (D)=\delta (K-D)$ and that $\beta (D)=\beta (K-D).$
The next lemma will allow us to study special divisors 
of degree $\leq g-1$. 

\begin{lem}
\label{residuel}
Let $D\in \Div(X)$ be an effective and special divisor of degree $d$.
\begin{description}
\item[({\cal{i}})] If
$$\dim\, |D |\leq \frac{1}{2} (d-\delta (D))+k$$
for a positive integer $k$, then
$$\dim\, |K-D |\leq \frac{1}{2} (\deg(K-D)-\delta(K-D))+k .$$
\item[({\cal{ii}})]
If $$\dim\, |D |\leq \frac{1}{2} (d-\beta (D))-k$$
for a positive integer $k$, then
$$\dim\, |K-D |\leq \frac{1}{2} (\deg(K-D)-\beta (K-D))-k .$$
\end{description}
\end{lem}

\begin{proof} It is a straightforward calculation using Riemann-Roch.
\end{proof}

The following lemma concerns non-trivial 
covering maps of degree $2$ between real curves.
\begin{lem}
\label{degre2}
Let $\varphi :X\rightarrow X'$ be a non-trivial covering map of degree $2$
between two real curves $X$ and $X'$. If there exists a real point 
$P\in X'(\RR )$
such that $\varphi^{-1} (P)=\{ P_1 , P_2\}$, with $P_1 $ and $P_2 $
real points not contained in the same connected component of $X(\RR)$,
then $\varphi (C_1 )=\varphi (C_2 )=C$
and $\varphi^{-1} (C)=C_1 \cup C_2$, with $C,C_1,C_2 $ the real
connected components containing the points
$P, P_1 , P_2$ respectively
\end{lem}

\begin{proof}
Since $C_1$ and $C_2$ are connected, we have $\varphi (C_1 )\subseteq C$
and $\varphi (C_2 )\subseteq C$. Moreover $\varphi (C_1 )$ and $\varphi (C_2 )$
are closed connected subsets of $C$ since $\varphi$ is proper. The morphism
$\varphi $ is étale at $P$ , hence there is an open neighbouroud $U$ of $P$
such that for any $Q\in U$ we have $\varphi^{-1} (Q)=\{ Q_1 ,Q_2 \}$
with $Q_i \in C_i $ for $i=1,2$. In fact, this situation does not change
when we run along $C$ since $C_1 \cap C_2 =\emptyset $ and thus $C$ 
cannot have a branch point.
\end{proof}

We state the main result of the paper.
\begin{thm}
\label{cliffreel}
Let $D$ be an effective and special divisor of degree $d$, and let
$k\in\NN$. Then either
$$\dim\, |D|\leq \frac{1}{2} (d-\delta (D))+k$$
or
$$\dim\, |D|\leq \frac{1}{2} (d-\beta (D))-k.$$
\end{thm}

\begin{proof}
Let $r=\dim\, |D|$. If we can show that
\begin{equation}
\label{equ1}
r\geq\frac{1}{2} (d-\delta (D))+k+1
\end{equation}
and
\begin{equation}
\label{equ2}
r>\frac{1}{2} (d-\beta (D))-k
\end{equation}
do not hold simultaneously, we shall have proved
the theorem. Assuming they are true, then we use Lemmas 
\ref{reductbasepointfree} and \ref{residuel} to restrict the
investigation for the case $D$ base point free and $0\leq d \leq g-1$.

By Lemma \ref{huisman2}, we have
\begin{equation}
\label{equ3}
d+\delta (D)\geq 2s+2k.
\end{equation}
Using (\ref{equ3}) and (\ref{equ2})
we obtain 
\begin{equation}
\label{equ4}
r> \frac{1}{2} (d-s+\delta (D))-k \geq \frac{1}{2}s.
\end{equation}
By (\ref{equ1}) and (\ref{equ2}) 
$$\delta (D)\geq d-2r+2k+2 ,$$
and $$\beta (D)=s-\delta (D)\geq d-2r-2k+1 .$$
Hence $$s\geq 2d-4r+3$$
and using (\ref{equ4}) we obtain
\begin{equation}
\label{equ5}
r> \frac{1}{3}(d+1).
\end{equation}

There are two cases to be looked at:

({\cal{i}}) First, $D$ is simple. \\
In this case, $X$ is mapped birationally by $|D|$
onto a curve of degree $d$ in $\PP_{\RR}^r$. According to a well-known
formula of Castelnuovo \cite[p. 116]{ACGH} for the genus of a curve in
$\PP_{\RR}^r$,
we have
\begin{equation}
\label{equ6}
g\leq m(d-1 -\frac{1}{2}(m+1)(r-1))
\end{equation}
where $m=[\frac{d-1}{r-1} ]$. By Clifford's theorem and (\ref{equ5}),
$m=2$ or $m=3$, since $r\geq 2$.
If $m=2$ (resp. $m=3$), replacing in (\ref{equ6}) and using
(\ref{equ5}), we get $d>g$ (resp. $d> g-1$), contradicting the fact
that $D$ was supposed of degree $\leq g-1$..
\air

({\cal{i}}) Second, $D$ is not simple.\\
Consider the map 
$f:X\rightarrow \PP_{\RR}^r$
associated to $|D|$. Let $X'$ be the normalization of $f (X)$. Then 
the induced morphism $\varphi :X\rightarrow X'$ is a non-trivial covering
map of degree $t\geq 2$ and there is $D'\in \Div(X')$ such that
$|D'|$ is a $g_{\frac{d}{t}}^{r}$ and such that $D=\varphi^* (D')$.
\air

To finish the proof, we proceed in three steps:
\air
{\sc Step 1}: $D'$ is non-special and $t=2$.\\
If $D'$ were a special divisor on $X'$, 
then $2r\leq \frac{d}{t} <\frac{3r}{t}$ (by (\ref{equ5}) and
Clifford's theorem), contradicting $t\geq 2$.
Hence $D'$ is non-special and $r=\frac{d}{t}-g'$ by
Riemann-Roch, where $g'$ denotes the genus of $X'$.
Using (\ref{equ5}), we get $0\leq g'<\frac{3r-1}{t}-r$
and thus, $t=2$.
\air

{\sc Step 2}: $X'$ is an $M$-curve, $\delta(D)=2\delta(D')=2g'+2$
and $k=0$.\\
Since $\varphi :X\rightarrow X'$ is a non-trivial covering
map of degree $2$ and $D=\varphi^* (D')$, we have 
$\delta(D)\leq 2\delta(D')\leq 2g'+2$ by Lemma \ref{degre2} and
by Harnack's inequality.
Since $r=\frac{1}{2}(d-2g')$, using the inequality (\ref{equ1})
we obtain $2g'+2k+2\leq \delta(D)$. Thus, $\delta(D)=2\delta(D')=2g'+2$
and $k=0$.
\air

{\sc Step 3}: $s=2g'+2$.\\
Let $P_1 ,\ldots ,P_{g'+1}$ be real points of the support of $D'$ such that
$E=\varphi^* (E'=P_1 +\ldots +P_{g'+1})$ satisfies $\delta(E)=2g'+2$.
By Riemann-Roch, $\dim\,|E'|\geq 1$, hence $\dim\,|E|\geq 1$. Using 
Lemma \ref{1parcomp}, $s=2g'+2$.

Summing up, $r=\frac{1}{2}(d-2g')\leq \frac{1}{2}d
=\frac{1}{2}(d-\beta (D))-k$ and this contradicts (\ref{equ2}).
\end{proof}

In the previous theorem, the case $k=0$ gives Theorem A announced
in the introduction.
\begin{thm}
\label{cliffreelA}
Let $D$ be an effective and special divisor of degree $d$. Then either
$$\dim\, |D|\leq \frac{1}{2} (d-\delta (D))\,\, (\Clif \, 1)$$
or 
$$\dim\, |D|\leq \frac{1}{2} (d-\beta (D))\,\, (\Clif \,2)$$
\end{thm}

\begin{rema}
{\rm The result of Theorem \ref{cliffreelA} is natural since,
for real curves without real points,
the inequalities $(\Clif\, 1)$ and $(\Clif\, 2)$ both become the classical
Clifford inequality. Moreover,
Theorem \ref{cliffreel} shows that in case $r=\dim\, |D|$ exceeds
the right hand term of $(\Clif\, 1)$ by $k>0$, 
then $r$ is exceeded by the right hand term of $(\Clif\, 2)$ by at least
$k-1$. It suggests that the inequalities $(\Clif\, 1)$ and $(\Clif\, 2)$ 
are not completely independant.}
\end{rema}

Classically, in the theory of special divisors, 
if the Clifford inequality becomes an equality for a divisor
different from $0$ and from the canonical divisor, then the curve is
hyperelliptic. We will show below that, for certain divisors on real
hyperelliptic curves, the inequalities
of Theorem \ref{cliffreelA} become equalities.
A real hyperelliptic curve is a real curve $X$ such that $X_{\CC}$
is hyperelliptic, i.e. $X_{\CC}$ has a $g_2^1$ 
(a linear system of dimension 
$1$ and degree $2$). Since this $g_2^1$ is unique, it is
a real linear system i.e. $X$ has a $g_2^1$ (see \cite[Lem. 4.2]{Mo}).
As always, we assume that $X(\RR )\not=\emptyset$ 
and moreover that $g\geq 2$.

\begin{prop} 
\label{hyper}
Let $X$ be a real hyperelliptic curve and let $D$ be an effective special
divisor of degree $d$ on $X$.
\begin{description}
\item[({\cal{i}})] If $\delta( g_2^1 )=0$ then
$$\dim\, |D|\leq \frac{1}{2} (d-\delta (D))\,\,\,\, (\Clif \, 1).$$
\item[({\cal{ii}})] If $\delta( g_2^1 )=2$ then $s=2$. Moreover either
$$\dim\, |D|\leq \frac{1}{2} (d-\delta (D))\,\,\,\, (\Clif \, 1)$$
or
$$\dim\, |D|= \frac{1}{2} (d-\beta (D))\,\,\,\, (\Clif \,2).$$
\end{description}
\end{prop}

\begin{proof}
By Lemma \ref{residuel}, we may assume that $d\leq g-1$.
Set $r=\dim\, |D|$.

Firstly, we consider that $\delta( g_2^1 )=0$. 
A consequence of the geometric version 
of the Riemann-Roch Theorem is that any complete and special 
$g_d^r$ on $X_{\CC }$
is of the form
$$rg_2^1+D',$$ where $D'$ is an effective divisor of degree $d-2r$
which has no fixed part under the hyperelliptic involution 
$\imath$ induced by the $g_2^1$. Since $\delta( g_2^1 )=0$,
we get $\delta( D )\leq \deg(D')=d-2r$. Hence ({\cal{i}})
of the proposition. Moreover, (Clif 1) becomes an equality if
and only if $D'=\sum_{i=1}^{\delta (D)} P_i$ with one $P_i$ in each component
of $X(\RR )$ where the
degree of the restriction of $D$ is odd.

Secondly, we assume that $\delta( g_2^1 )=2$. By Lemma \ref{degre2}
or the proof of {\sc Step 3} of the previous theorem, $s=2$ and the 
hyperelliptic involution exchanges the two connected
components of $X(\RR )$. If $r$ is even, the proof runs as in the case
$\delta( g_2^1 )=0$ and we get the inequality (Clif 1).
If $r$ is odd, we again write $$|D|=rg_2^1+D',$$ 
where $D'$ is an effective divisor of degree $d-2r$
which has no fixed part under the hyperelliptic involution.
If $\deg(D')\geq 2$ or if $\delta(D')\geq 1 $, then the inequality 
(Clif 1) works. If not, $D$ is the inverse image, 
by the morphism associated to the $g_2^1$,
of an effective divisor
$D''$ on $\PP_{\RR}^1$.
Consequently $|D|=rg_2^1$ and 
$$\dim\, |D|= \frac{1}{2} (d-\beta (D)).$$
\end{proof}

\begin{rema}
{\rm A real hyperelliptic curve such that $\delta( g_2^1 )=2$ is given
  by the real polynomial equation $y^2=f(x)$, where $f$ is a monic
  polynomial of degree $2g+2$, with $g$ odd, and where $f$ has no real
  roots \cite[Prop. 4.3]{Mo}.}
\end{rema}

Let $D$ be a special and effective divisor of degree $d$ 
on a real curve $X$ such that $\dim\, |D|= \frac{1}{2} (d-\beta (D))-k$
with $k\in\NN$. According to Theorem \ref{cliffreel}, we have 
$\dim\, |D|\leq \frac{1}{2} (d-\delta (D))+k+1$. Similarly,
if $\dim\, |D|= \frac{1}{2} (d-\delta (D))+k+1$ with $k\in\NN$ then, 
using  Theorem \ref{cliffreel},
$\dim\, |D|\leq \frac{1}{2} (d-\beta (D))-k $. 
Consequently, we will say that $D$ is extremal (for
the real Clifford inequalities) if $\dim\, |D|= \frac{1}{2} (d-\beta
(D))-k= \frac{1}{2} (d-\delta (D))+k+1$ for some $k\in\NN$. 

\air
Looking at the example of hyperelliptic curves, we may ask the following 
questions:
\air

Do there exist extremal divisors and what geometric
properties does this imply for $X$?
\air

In case $k=0$ and $D$ is extremal, does it follow that
$X$ is an hyperelliptic curve
with $\delta( g_2^1 )=2$?
\air

Before giving an answer to these questions in Theorem \ref{siegalite},
we state some classical
results concerning extremal complex curves and special divisors on
complex curves that easily extend to real curves.
Recall that a non-degenerate curve $X$ in $\PP_{\RR}^r$ is called
extremal if the genus is maximal with respect to the degree of $X$
(cf. \cite[p. 117]{ACGH}).

\begin{lem} (\cite[Lem. 3.1]{ELMS})
\label{accola}
Let $D$ and $E$ be divisors of degree $d$ and $e$ on a curve $X$ of
genus $g$ and suppose that $|E|$ is base point free. Then 
$$\ell (D)-\ell (D-E)\leq \frac{e}{2}$$
if $2D-E$ is special.
\end{lem}

The lemma applies in case $D$ is semi-canonical i.e. $2D=K$.

\begin{lem} (\cite[Lem. 1.2.3]{CKM})
\label{multiple}
Let $g_k^1$ be a base point free pencil on a curve $X$ of genus
$g$. Assume
$k(k-1)\leq 2g-2$ and let $n\in \NN$ such that $n\leq
\frac{2g-2}{k(k-1)}$.
Then $\dim\, ng_k^1=n$ and $ng_k^1$ is base point free.
\end{lem}

\begin{lem} (\cite{Ac}, \cite{Co-Ma} p. 200 and \cite{ACGH} p. 122)
\label{extremale} 
Let $X$ be an extremal curve of degree $d>2r$ in $\PP_{\RR}^r$ $(r\geq
3)$. Then one of the followings holds:
\begin{description}
\item[({\cal{i}})] $X$ lies on a rational normal scroll $Y$ in
  $\PP_{\RR}^r$ ($Y$ is real, see \cite{ACGH} p. 120). Write
  $d=m(r-1)+1+\varepsilon$ where
  $m=[\frac{d-1}{r-1} ]$ and
  $\varepsilon\in \{0,2,\ldots ,r-2\}$. $X_{\CC}$ has only finitely many
  base point free pencils of degree $m+1$ (in fact, only $1$ for
  $r>3$, and $1$ or $2$ if $r=3$). These pencils are swept out by the
  rulings of $Y_{\CC}$. Moreover $X_{\CC}$ has no $g_m^1$.
\item[({\cal{ii}})] $X$ is the image of a smooth plane curve $X'$ of
  degree $\frac{d}{2}$ under the Veronese map
  $\PP_{\RR}^2\rightarrow\PP_{\RR}^5$.
\end{description}
\end{lem}

We give an answer to the questions of the previous page 
in the following theorem.
\begin{thm}
\label{siegalite}
Assume $D$ is an effective and special divisor of degree $d$ such that
$$r=\dim\, |D|= \frac{1}{2} (d-\delta (D))+k+1= \frac{1}{2} (d-\beta (D))-k$$
i.e. $D$ is extremal.
Then $k=0$, $X$ is an hyperelliptic curve with 
$\delta( g_2^1 )=2$ and $|D|=rg_2^1 $ with $r$ odd.
\end{thm}

\begin{proof}
We set $\delta=\delta (D)$ and $\beta =\beta (D)$.

We now copy the proof of Theorem \ref{cliffreel}, as in that situation
we may assume $D$ is base point free (Lemma \ref{reductbasepointfree})
and $0\leq d \leq g-1$ (Lemma \ref{residuel}).

The inequality
(\ref{equ3}) remains valid. The inequalities (\ref{equ1}) and (\ref{equ2})
are now equalities
\begin{equation}
\label{equ7}
r=\frac{1}{2} (d-\delta )+k+1.
\end{equation}
and
\begin{equation}
\label{equ8}
r=\frac{1}{2} (d-\beta )-k.
\end{equation}
By (\ref{equ8}), $s$ is even.
The inequality (\ref{equ4}) becomes
$r\geq \frac{1}{2} s$.
Using (\ref{equ7}) and (\ref{equ8}), we obtain
$\delta = d-2r+2k+2 ,$ and
$\beta =s-\delta = d-2r-2k $.
Hence $s=2d-4r+2$.
Using (\ref{equ4}) again, we obtain
\begin{equation}
\label{equ9}
r\geq \frac{1}{3}(d+1).
\end{equation}

We have one of the two following
possibilities:
\air

({\cal{i}}) $D$ is simple. \\
The linear system $|D |$ embeds $X$ in $\PP_{\RR}^r$
as a curve of degree $d$. Using the facts that $d\leq g-1$ and 
$2r\leq d\leq 3r-1$ (by (\ref{equ9}) and Clifford's theorem)
a straightforward calculation shows that the Castelnuovo's inequality
(\ref{equ6}) is an equality and that $m\geq 3$, $d=3r-1$ and $g=3r$ is the only
possibility.
By (\ref{equ8}) and (\ref{equ7}), we have $\delta =r+1+2k$ and 
$s-\delta =r-1-2k$. Hence $s=2r$. By \cite[Lem. 2.9]{Ac2}, $D$ is
semi-canonical i.e. $2D=K$.
\air

At this moment of the proof there is no contradiction about the
existence of such extremal and simple divisor $D$. The geometric
properties of extremal curves will give this contradiction.
\air

{\sc Case 1}: $r=2$.\\
We identify $X$ via $|D|$ with a smooth plane quintic curve. The
contradiction is given by 
$\delta\geq 3$ and the fact that $X$ has a unique
pseudo-line (the definition of a pseudo-line is in the next section).
\air

{\sc Case 2}: Either $r> 5$ or  $r=5$ and $X$ is not a smooth plane curve.\\
Then $X_{\CC}$ has a unique $g_4^1$ (Lemma \ref{extremale}). 
Hence this $g_4^1$ is real
(cf. \cite{Mo}) i.e. there is an effective divisor $E$ of degree $4$
such that $|E |=g_4^1 $. By Lemma \ref{multiple}, $2E$ is base point
free and $\dim\, |2E|=2$. Let $r'=\dim\, |D-2E|$. Since $D$ is
semi-canonical, we have $r'\geq 5-4=1$ by Lemma \ref{accola}.
Hence $D-2E$ is special, moreover $\delta (D-2E)=\delta$ and $D-2E$ is
also extremal for the same $k\in \NN$. Since $X$ is not hyperelliptic,
$D-2E$ is simple (see the part of the proof concerning non-simple
extremal divisors).
Hence $g=3r'$ and we get a
contradiction, since $r'<r$.
\air

{\sc Case 3}: $r=5$ and $X$ is a smooth plane curve.\\
By Lemma \ref{extremale}, $X$ is the image of a smooth plane curve of
degree $7$ under the Veronese embedding
$\PP_{\RR}^2\rightarrow\PP_{\RR}^5$.
Hence $X$ has a unique very ample $g_7^2 =|E|$. Using Lemma
\ref{accola} ($D$ is semi-canonical) and since 
$E$ calculates the Clifford index of $X_{\CC}$
(cf. Section 4 for the definition of the Clifford index), we have
$|D|=|2E|$. So $\delta =0$, which is impossible.
\air

{\sc Case 4}: $r=4$.\\
Similarly to {\sc Case 2}, $X$ has a $g_4^1 =|E|$. Let
$D'=D-E$. Applying Lemma \ref{accola}, we get $\dim\, |D' |\geq 2$.
By Riemann-Roch $\ell (K-(2D'-E))=\ell (2D'-E)-10+12-1 >0$.
Consequently, according to  Lemma \ref{accola}, $\ell (D'-E)=\ell
(D-2E)>0$. More precisely either $\ell
(D-2E)=1$ or $X$ would have a $g_3^1$ contradicting Lemma
\ref{extremale}. So $|D|=|2E+D''|$, with $D''$ an
effective divisor of degree $3$. Hence $\delta \leq 3$, which is 
again impossible.
\air

{\sc Case 5}: $r=3$ and $X_{\CC}$ has a unique $g_4^1 =|E|$.\\
By Lemma \ref{accola} and since $X$ has no $g_3^1$, we get
$|D|=|2E|$ and a contradiction on $\delta$.
\air

{\sc Case 6}: $r=3$ and $X_{\CC}$ has two $g_4^1$, $|E|$ and $|F|$
which are real.\\ 
We know that $X$ lies on a unique quadric $S$, and
$|E|$ and $|F|$ correspond to the rulings of $S$. More precisely $X$ is
of bi-degree $(4,4)$ on $S\cong
\PP_{\RR}^1 \times \PP_{\RR}^1$. Consequently $\delta (E)=\delta(F)$.
By Lemma \ref{accola} for $D$ and $E$, and since $X$ has no $g_3^1$,
we get $|D|=|2E|$ or $|D|=|E+F|$.
Hence $\delta =0$, contradiction.
\air

{\sc Case 7}: $r=3$ and $X_{\CC}$ has two $g_4^1$, $|E|$ and $|\bar{E}|$
which are complex and switched by the complex conjugation.\\
We argue similarly as in the previous case, but on the complex curve
$X_{\CC}$. We obtain that $|D|=|E+\bar{E}|$ i.e. $\delta =0$, which is 
impossible.
\air

({\cal{ii}}) $D$ is not simple.\\
Here $|D|$ induces
a non-trivial covering
map $\varphi :X\rightarrow X'$ of degree $2$
on an $M-$curve $X'$
of genus $g'$ (see the proof of Theorem \ref{cliffreel}). 
There is an effective divisor $D'\in \Div(X')$ such that
$|D'|$ is a $g_{\frac{d}{2}}^{r}$ and such that $D=\varphi^* (D')$.
Moreover, following the proof of Theorem \ref{cliffreel}, 
we see that $D'$ is non-special, $r=\frac{1}{2}(d-2g' )$,
$\delta =2\delta(D')=2g'+2 =s$. The identities
(\ref{equ7} ) and (\ref{equ8} ) say that $k=0$ and that $g'=0$
i.e. that $X$ is an
hyperelliptic curve. By Proposition \ref{hyper} we get that $\delta(
g_2^1 )=s=2$ and $|D|=rg_2^1 $ with $r$ odd (if $r$ is even, it
contradicts (\ref{equ7})).
\end{proof}

The Clifford type inequalities from Theorem \ref{cliffreel} seem to be
the best possible since in the previous proof, the extremal cases for
these inequalities correspond to extremal Castelnuovo curves. As in
the complex situation, these inequalities become equalities in
non-trivial cases, only if the curves are hyperelliptic.

From Theorem \ref{cliffreel} and Theorem \ref{siegalite}, 
we may derive the following result which corresponds to Theorem B in
the introduction.

\begin{thm}
\label{meilleure}
Let $X$ be a real curve.
Let $D\in \Div(X)$ be an effective and special divisor of degree $d$.
\begin{description}
\item[({\cal{i}})] If $X$ is hyperelliptic then, either
$$\dim\, |D|\leq \frac{1}{2} (d-\delta (D))\,\,\,\, (\Clif \, 1)$$
or 
$$\dim\, |D|\leq \frac{1}{2} (d-\frac{1}{2}(s-2)),$$
\item[({\cal{ii}})] else, either
$$\dim\, |D|\leq \frac{1}{2} (d-\delta (D))\,\,\,\, (\Clif \, 1)$$
or
$$\dim\, |D|\leq \frac{1}{2} (d-\frac{1}{2}(s-1)).$$
\end{description}
\end{thm}

\begin{proof}
Let $r=\dim\, |D|$. Assume $r=\frac{1}{2} (d-\delta (D))+k+1$. By
Theorem \ref{cliffreel} and Theorem \ref{siegalite}
we get $r\leq \frac{1}{2} (d-\beta (D))-k$ and the inequality is
strict if $X$ is not hyperelliptic.
Hence $2k+1\leq \frac{1}{2}(\delta(D)-\beta(D))=\delta(D)-\frac{1}{2}
s$ with equality only if $X$ is hyperelliptic, completing the proof.
\end{proof}

We show now that the inequalities of Theorem \ref{meilleure} may
become equalities.
\begin{exem}
{\rm Let $X$ be  an hyperelliptic curve $X$ such that $\delta (g_2^1
  )=2$. If $D$ is an element of the $g_2^1$, then $D$ does not
  satisfy the inequality ($\Clif \, 1$) and the second inequality of
  Theorem \ref{meilleure} ({\cal{i}}) is an equality.

Let $X$ be a real trigonal curve, i.e. $X$ has a $g_3^1$. We assume
that $\delta (g_3^1 )=3$ and we take $D$ an element of the $g_3^1$. By
\cite[p. 179]{G-H}, such a trigonal curve exists. Then $D$ does not
satisfy the inequality ($\Clif \, 1$), but it gives an example of a
divisor for which equality holds in the second inequality of
Theorem \ref{meilleure} ({\cal{ii}}).}
\end{exem}

\section{Special real curves in projective spaces}

Let $X\subseteq \PP_{\RR}^{r}$, $r\geq 2$, 
be a smooth real curve, $X$ is
non-degenerate if $X$ is not contained  in an hyperplane of
$\PP_{\RR}^n$. We assume, in what follows, that $X$ is non-degenerate.
We say that $X$ is special (resp. non-special) if the divisor
associated to the sheaf of hyperplane sections $\L_{X} (1)$
is special (resp. non-special).

Let $C$ be a connected component of $X(\RR )$.
The component $C$ is called a pseudo-line if the canonical class of $C$
is non-trivial in $H_1 (\PP_{\RR}^r (\RR ),\ZZ /2)$.
Equivalently, $C$ is
a pseudo-line
if and only if for each real hyperplane $H$,
$H(\RR )$ intersects $C$ in an odd number of points, 
when counted with multiplicities (see \cite{Hu}).

In this section, we wish to discuss some conditions 
under which we may bound the genus, the number of pseudo-lines,
and the number of ovals of a non-degenerate 
smooth real curve in $\PP_{\RR}^{r}$.
For the genus, if $X$ is a smooth plane curve of degree $d$, we have 
$$g=\frac{1}{2} (d-1)(d-2).$$
When $r\geq 3$, there is no formula for the  genus of $X$ 
in terms of its degree.
The situation is therefore more complicated. However, there is an inequality
of Castelnuovo (inequality (\ref{equ6})) that we have already seen
in the proof of Theorem \ref{cliffreel}. 

The following proposition improves the Castelnuovo inequality
for non-special real curves.
\begin{prop}
\label{nonspecial}
Let $r\geq 2$ be an integer and $X\subseteq \PP_{\RR}^r$ be a non-degenerate
real curve. Let $d$ be the degree of $X$ and
$\delta$ (resp. $\beta$) be the number of pseudo-lines 
(resp. ovals) of $X$. Assume $d+2k<2r+\delta$ and $d-2k<2r+\beta$ for
some $k\in \NN$. Then
$$g\leq d-r,$$ and equality holds if and only if $X$ is linearly normal
i.e. if and only if the restriction map
$$H^0 (\PP_{\RR}^r ,\O (1))\rightarrow H^0 (X ,\O_X (1))$$
is surjective.
\end{prop}

\begin{proof}
Let $H$ be a hyperplane section of $X$ i.e. a divisor obtained by cutting
out the curve by a real hyperplane. Then 
$\dim\,|H|\geq r>\frac{1}{2}(d-\delta(H))+k$ and 
$\dim\,|H|\geq r>\frac{1}{2}(d-\beta(H))-k$
by the hypotheses.
Theorem \ref{cliffreel} says that $H$ is non-special and by Riemann-Roch, 
$$g=d-\dim\,|H|\leq d-r.$$
Clearly, the previous inequality becomes an equality if and only
if the map $H^0 (\PP_{\RR }^r ,\O (1))\hookrightarrow
H^0 (X ,\O_{X} (1))$ is an isomorphism.
\end{proof}

One may wonder what can be said about the number of pseudo-lines and
ovals of $X$ when $X\subseteq \PP_{\RR }^r$ is a non-special real curve.
The following proposition shows that there is no restriction on these numbers
except the fact that the number of pseudo-lines should be congruent to
the degree of the curve modulo 2.

\begin{prop}
\label{nonspecialinproj}
Let $r\geq 3$ be an integer and $X$ be a 
real curve. Let $\delta $ be an integer $\leq s$.
There is a smooth embedding $\varphi :X\hookrightarrow \PP_{\RR}^r$
such that $X$ is non-special in $\PP_{\RR}^r$
and $X$ has $\delta$ pseudo-lines provided that
$\delta =g+r\, \mod\, 2$.
\end{prop}

\begin{proof}
Since $\delta < g+r$ and $\delta =g+r\, \mod\, 2$, 
there is an effective divisor $D$ of degree
$g+r$ such that $\delta (D)=\delta$.
Choosing $D$ general, $D$ is non-special
and $D$ is very ample 
(see the proof of a theorem of Halphen \cite[p. 350] {Ha}).
The morphism associated
to $|D|$ gives the result.
\end{proof}

We show now that Theorem \ref{cliffreelA} gives a lower bound 
on the number of ovals
or the number of pseudo-lines of special real space curves.
\begin{prop}
Let $X\subseteq \PP_{\RR}^r$ be a special non-degenerate real curve
of degree $d$.
Let $\delta$ (resp. $\beta$) denote the number of pseudo-lines
(resp. ovals) of $X$. Then, either
$$\beta \geq r+s-g+1$$
or 
$$\delta \geq r+s-g+1 .$$
\end{prop}

\begin{proof}
By Theorem \ref{cliffreelA}, we have two possibilities since the hyperplane
section of $X$ is special.

Firstly, $d\geq 2r +\delta =2r +s -\beta $. Hence $\beta \geq (r+s)+(r-d)$.
But $g>d-r$ since $X$ is special, so $\beta \geq r+s -g+1$.

Secondly, $d\geq 2r +\beta$ and by a similar argument we get
$\delta \geq r+s -g+1$.
\end{proof}

In particular, for $(M-2)$-curves, the above proposition gives:
\begin{cor}
Let $X\subseteq \PP_{\RR}^r$ be a special non-degenerate $(M-2)$-curve
of degree $d$.
Let $\delta$ (resp. $\beta$) denote the number of pseudo-lines
(resp. ovals) of $X$. Then, either
$\beta \geq r$
or 
$\delta \geq r.$
\end{cor}

\section{Distribution of the special divisors between the inequalities
(Clif 1) and (Clif 2)}

The special divisors splits naturally into two a priori equivalent sets:

1) divisors $D$ with $\delta (D)\leq \frac{s}{2}$: they satisfy (Clif 1).

2) divisors $D$ with $\delta (D)> \frac{s}{2}$: they satisfy (Clif 2).

At this level the inequalities seem to be equivalent. In fact,
according to the results of section 2 for hyperelliptic curves, the inequality 
(Clif 1) seems to be dominant. Moreover, every real 
curve with non-empty real part has a special divisor 
satisfying (Clif 1) but not (Clif 2): the canonical divisor. 
We will show that there are real curves
for which every special divisor satisfies (Clif 1).

\subsection{ Clifford inequalities and the number of real connected components}

We state the following Huisman's result \cite[Th. 3.2]{Hu}.

\begin{thm} Assume $X$ is an $M$-curve or an $(M-1)$-curve.
Let $D\in \Div(X)$ be an effective and special divisor of degree $d$.
Then $$\dim\, |D|\leq \frac{1}{2} (d-\delta (D))\,\,\,\, (\Clif \,1).$$
\end{thm}

Unfortunately, we have seen that the inequality (Clif 1) is not 
valid for all real curves, but the previous theorem also shows that (Clif 1)
is dominating (Clif 2).

The following theorem gives Huisman's theorem for $n=0$ and $n=1$.
\begin{thm}
\label{huisman3}
Let $X$ be an $(M-n)$-curve with $n$ an integer such that
$0\leq n\leq g+1$. If $D$ is an 
effective and special divisor of degree $d$ on $X$ then
$$\dim\, |D|\leq \frac{1}{2} (d-\delta (D))+[\frac{n}{2}] .$$
\end{thm}

\begin{proof}
Let $D$ be an effective and special divisor of degree $d$ such that
$\dim\,|D|=\frac{1}{2} (d-\delta (D))+k+1$.
By Lemma \ref{residuel}, we may assume that $d\leq g-1$.
From Lemma \ref{huisman2}, we have $d+\delta(D)\geq 2s+2k$.
Hence $g-1\geq d \geq 2s-\delta (D)+2k\geq s +2k$.
Consequently $$k\leq \frac{1}{2}(g-1-s)$$
and the proof is done.
\end{proof}

\subsection{Clifford inequalities and dimension of linear systems}

The following proposition shows that for linear systems of 
dimension $r$ which do not satisfy (Clif 1), the excess
can be bounded in terms of $r$.

\begin{prop}
Let $D$ be an effective and special divisor of degree $d$ on a
real curve $X$. Assume that $r=\dim\, |D|= \frac{1}{2} (d-\delta (D))+k+1.$
Then $k\leq [\frac{r-1}{2}]$.
\end{prop}

\begin{proof}
Assume $r=2n-\varepsilon$ with $\varepsilon\in\{ 0,1\}$.
We proceed by induction on $n$.

If $r=1$, we get $2=d-\delta (D)+2k+2$. Since $d\geq \delta (D)$, 
we have $k=0$ and $d=\delta (D)$.

If $r=2$, we get $4=d-\delta (D)+2k+2$ i.e. $2=d-\delta (D)+2k$.
If $k\geq 1$, we must have $d=\delta (D)$, hence $d+\delta (D)\leq 2s$
but Lemma \ref{huisman2} says that $k=0$, which is impossible.
So $k=0$ and $d=\delta (D)+2$.

Assume $n>1$ and $k>0$. Choose two real points $P_1 ,P_2$ in the same real 
connected component such that $\dim\, |D-P_1 -P_2 |=r-2$.
Then $r-2=2(n-1)-\varepsilon =\frac{1}{2} ((d-2)-\delta (D))+(k-1)+1.$
By the induction hypothesis, $(k-1)\leq [\frac{(r-2)-1}{2}]$ i.e.
$$k\leq [\frac{r-1}{2}].$$
\end{proof}

\begin{cor}
A pencil (i.e. a complete $g_d^1$) which does not satisfy (Clif 1)
is composed by divisors of the form $P_1 +\ldots +P_s$
with $P_1 ,\ldots ,P_s $ some real points of $X$ 
such that no two of them belong to the same connected component
of $X(\RR )$.
\end{cor}

\begin{proof}
Looking at the proof of the previous proposition, we see that
$d=\delta(D)$. Using Lemma \ref{1parcomp}, we get the stated result.
\end{proof}

\subsection{ Clifford index of real curves}

In this paragraph, we will introduce the notion of the
Clifford index of a real curve. The Clifford index of 
a complex curve is a classical concept in the theory of complex curves. So
we adapt the definitions in the real case.

Let $D$ be an effective and special divisor of degree $d$ on $X$ 
such that $\dim\, |D|=r$. The expression $d-2r$ is called
the Clifford index of $D$ and is denoted by $\Cliff (D)$.
For a given degree,
the smaller the Clifford index of $D$ is, the more global sections 
$\L (D)$ has. We restrict attention to divisors with both
$\ell (D)>1$ and $\ell (K-D)>1$. We say that these divisors contribute to 
the Clifford index of $X$ (denoted by  $\Cliff(X)$), 
which is defined as the minimum of their 
Clifford indices. If such divisors do not exist, we will
say that the Clifford index of $X$ is infinite.
We say that $D$ computes the Clifford index of $X$
if $D$ contributes to the Clifford index of $X$ and 
$\Cliff (D)=\Cliff(X)$.
The notion of Clifford index of $X$ refines the notion of 
gonality of $X$, which is the smallest degree of a map from $X$ to 
$\PP_{\RR}^1$.

We use the same definitions
concerning the notions of Clifford index and gonality
when $Y$ is a smooth connected projective curve over the complex numbers and
$D$ is a divisor on $Y$.
From Brill-Noether theory, it follows that a complex curve
has Clifford index $\leq [\frac{(g-1)}{2} ]$ and gonality 
$\leq [\frac{(g+3)}{2} ]$ with equalities if the curve is general.
\air

What about real curves ?
\air

Of course, if the Clifford index of $X$ is finite then
$\Cliff (X)\geq \Cliff (X_{\CC})$.
By \cite[Th. 5.4]{Mo}, we know that a real curve has always a $g_g^1 $ and we
conjecturate that there exists a $g$-gonal real curve of genus $g$
for any $g\geq 2$. 
The conjecture is proved for $g\leq 4$ \cite[Prop. 5.11]{Mo}.

\begin{prop} 
\label{nocliff}
A $g$-gonal real curve has an infinite Clifford index.
\end{prop}

\begin{proof}
The proof is clear since the residual linear system of a $g_g^1 $
is zero dimensional. Consequently, a $g$-gonal real curve 
does not have any divisor contributing to the Clifford index.
\end{proof}

Concerning the existence of the Clifford index of a real curve,
we state the following results:
\begin{prop}
\label{existence1}
Let $X$ be a real curve such that $X_{\CC }$ is $k$-gonal
with $k\leq \frac{1}{2} (g+1)$. Then the Clifford 
index of $X$ is $\leq g-3$.
\end{prop}

\begin{proof}
By \cite[Th. 5.4]{Mo} $X$ has gonality $\leq 2k-2\leq g-1$. Hence
there exists an effective divisor $D$ on $X$ of degree $g-1$
such that $\dim\, |D|=1$. Since $\dim\, |K-D|=1$, $D$ contributes to 
the Clifford index of $X$, and $\Cliff (X)\leq \Cliff(D)=g-3$.
\end{proof}

\begin{prop}
\label{existence2}
Let $X$ be a real curve such that $X_{\CC }$
has a unique linear system calculating $\Cliff (X_{\CC })$.
Then $\Cliff (X)=\Cliff (X_{\CC })$.
\end{prop}

\begin{proof} Let $c=\Cliff (X_{\CC })$ and let $D$ be an effective divisor 
on $X_{\CC}$ such that $\Cliff(D)=c$. Assume $\dim\, |D|=r$
and let $P_1 ,\ldots , P_r$ be real points (seen as closed points of $X_{\CC}$)
such that $\dim\, |D-\sum_{i=1}^{r} P_i |=0$. Hence $D-\sum_{i=1}^{r} P_i $
is linearly equivalent on $X_{\CC }$ to a unique effective divisor $D'$.
Since $X$ is a real curve, we have $\dim\, |\sum_{i=1}^{r} P_i
+\overline{D'}|=r$, so $\sum_{i=1}^{r} P_i + \bar{D'}$ is calculating
$\Cliff (X_{\CC })$. By the uniqueness property of $D$, it follows that 
$\sum_{i=1}^{r} P_i +\overline{D'}\in |\sum_{i=1}^{r} P_i +D'|$,
and in consequence $D'$ is linearly equivalent to $\overline{D'}$.
Since $\dim\, |D'|=0$, we have $D'=\overline{D'}$ and $D$
can be represented by a real divisor. The proof is now straightforward.
\end{proof}

One application of the previous proposition concerns the Clifford
index of a real curve whose complexification is a general $k$-gonal
curve.

\begin{cor}
\label{general}
Let $X$ be a real curve such that $X_{\CC }$ is a general $k$-gonal curve
with $3 \leq k< \frac{1}{2} g$. Then 
$$\Cliff (X)=\Cliff (X_{\CC })=k-2.$$
\end{cor}

\begin{proof}
From \cite[Example 3.3.4]{Co-Ma}, we see that on $X_{\CC }$ there is only one
linear system computing the Clifford index: the unique $g_k^1$.
By Proposition \ref{existence2}, this $g_k^1$ is real
and it calculates the Clifford index of $X$.
\end{proof}

In the following propositions, we investigate the relations between
the Clifford inequality (Clif 1) and the Clifford index of a real curve.
\begin{prop}
\label{index1}
Let $X$ be a real curve such that $X$ has an infinite Clifford index or
$\Cliff (X)\geq s$. 
Let $D\in \Div(X)$ be an effective divisor of degree $d$
then $$\dim\, |D|\leq \,\frac{1}{2} (d-\delta (D))\,\,\, (\Clif\, 1).$$
\end{prop}

\begin{proof}
Let $D$ be an effective and special divisor
which does not contribute to the Clifford index of $X$.
It means that $\dim\, |D|=0$ or $\dim\, |K-D|=0$.
Since $D$ and $K-D$ are effective ($D$ is special),
we have that either $D$ or $K-D$ satisfies (Clif 1). By
Lemma \ref{residuel}, $D$ verifies (Clif 1).

If $X$ has an infinite Clifford index, the proof follows from
what we have just said, since there is no special divisor contributing to
the Clifford index of $X$.

Assume now $\Cliff (X)\geq s$ and let $D$ be an effective and special divisor
of degree $d$ contributing to the Clifford index of $X$.
Then $\Cliff (D)=d-2\dim\, |D|\geq s$. So 
$\dim\, |D|\leq \,\frac{1}{2} (d-s)\leq \frac{1}{2} (d-\delta (D))$.
\end{proof}

\begin{cor}
\label{generale1}
Let $X$ be a real curve such that $X_{\CC }$ is general
and such that $s\leq [\frac{(g-1)}{2} ]$. 
Let $D\in \Div(X)$ be an effective divisor of degree $d$
then $$\dim\, |D|\leq \,\frac{1}{2} (d-\delta (D))\,\,\, (\Clif\, 1).$$
\end{cor}

\begin{prop}
\label{index2}
Let $X$ be a real curve such that $X$ has a finite Clifford index 
$c=\Cliff (X)$.
Let $D\in \Div(X)$ be an effective divisor of degree $d$
then
\begin{description}
\item[1)] $\dim\, |D|\leq \,\frac{1}{2} (d-\delta (D))$ if $\delta(D) \leq c$ 
and,
\item[2)] $\dim\, |D|\leq \,\frac{1}{2} (d-\delta (D))+
[\frac{1}{2}(\delta (D)-c)]$ if $\delta (D) > c$.
\end{description}
\end{prop}

\begin{proof}
Let $D$ be an effective and special divisor
of degree $d$. Looking at the proof of Proposition \ref{index1},
we may assume that $\delta (D) > c$ and that $D$ contributes
to the Clifford index of $X$. 

Assume $\dim\, |D|= \,\frac{1}{2} (d-\delta (D))+k$. Then
$\Cliff (D)=\delta (D)-2k\geq c$ and the proposition follows.
\end{proof}

For general curves, we complete Corollary \ref{generale1}
as follows:
\begin{cor}
\label{generale2}
Let $X$ be a real curve such that $X_{\CC }$ is a general curve
and such that $s > [\frac{(g-1)}{2} ]$. 
Let $D\in \Div(X)$ be an effective divisor of degree $d$
then $$\dim\, |D|\leq \,\frac{1}{2} (d-\delta (D))
+[\frac{1}{2}(s-\frac{(g-1)}{2})].$$
\end{cor}

\end{document}